\documentclass[12pt]{article}
\usepackage{amsfonts}
\usepackage{amssymb}
\usepackage{mathtools}
\usepackage{verbatim}
\usepackage{amsthm}
\theoremstyle{plain}
\newtheorem{theorem}{Theorem}
\newtheorem*{theorem*}{Theorem}

\newtheorem{lemma}{Lemma}
\newtheorem{definition}{Definition}
\newtheorem{definition*}{Definition}
\newtheorem{remark}{Remark}
\newtheorem{example}{Example}

\begin{document}

\begin{center}
{{\Large \bf On optimization on ravine functions.
 Minkowski-Cohn moduli surface in Cohn parameterization

}
}

\end{center}

\begin{center}
{\bf Nikolaj M. Glazunov } \end{center}

\begin{center}
{\rm Glushkov Institute of Cybernetics NASU, Kiev, } \\

{\rm  Email:} {\it glanm@yahoo.com }
\end{center} 


\bigskip


\begin{abstract}
This paper presents a brief overview of ravine functions using the example of the Minkowski-Cohn moduli surface from the point of view of optimization on it.
Elements of representation and solution of the minimization problem at a point are presented.
\end{abstract}

Ravine optimization methods were introduced by I.M. Gelfand and M. L. Tsetlin \cite{Co:MC}  in 1961. The author learned about these methods in the 1970s, while interning at the Computing Center of the USSR Academy of Sciences, Moscow State University, and the Steklov Mathematical Institute of the USSR Academy of Sciences, and participating in scientific seminars at these organizations.The first internship, lasting approximately a year and a half, took place in 1970–1971 during the author's postgraduate studies at the Institute of Cybernetics, while the second internship—also lasting about a year and a half—took place in 1974–1975, following the author's successful defense of their Candidate of Sciences dissertation.
   One of the properties of ravine functions is their extended level curves. We will consider the class of real ravine functions that are convex upward or downward. Upwardly convex functions with strongly elongated surfaces (and corresponding curves) of levels could be called mountainous or elevation functions, but we will also call such functions ravine functions.

Revine functions arise in the problems of number theory, approximation theory, cybernetics, discrete analysis, in problems related to the control processes of physiological mechanisms. 
The impressive optimization methods proposed by Gelfand and Tzeitlin, as well as the methods developed by their followers, refer to the optimization of a function at a point.
In this regard, we note the methods of transforming space  proposed by N.Z. Shor and his 
 student M.G. Zhurbenko \cite{shzh} and developed by N.Z. Shor  \cite{shor1,shor} 
and his students.
However, in a number of problems in number theory, approximation theory, cybernetics, discrete analysis, and applied research, the problem of finding a curve on the surface of a function
 that yields optimal values 
 of the function 
when one of the variables in the functions domain changes arises in the moduli spaces of the objects being studied. 
  In a number of cases, both the surface itself and the curves of minima or maxima on it may be non-differentiable.

This paper is devoted to the study and possible methods for solving such a problem, using the example of a function whose domain of value is a Minkowski-Cohn surface. Although the outline for constructing such a surface was given in Minkowski's well-known monograph \cite{Mi:DA}, and the surface itself in the $(\tau, \sigma)$ parameterization was constructed by Cohn \cite{Co:MC} in 1950, the study of this surface and problems on it continues to attract the attention of researchers. In particular, to the author's knowledge, the curve of maxima on this surface is unknown. Researchers from Austria and Australia recently (2025-2026) contacted the author regarding results and references to the author's work in this area, and the author provided such references.

The content of the paper is as follows. We present the definition of the function defining the Minkowski-Cohn surface and construct elements of optimization theory for the first  parameterization of this surface. The first parameterization is called the Cohn parameterization. The second and third parameterizations are called  Malyshev parameterizations and will be considered in the next paper. To find the curve of minima on this surface, it is convenient to use one or another parametrization depending on the subdomain of the Minkowski-Cohn surface under consideration. For known results see 
\cite{GM:MM,GM:P2,GGM:PM,gopms} and references therein.\\

\section{On non-differentiable functions and manifolds of ravine type}
    
   We will consider non-differentiable ravine-type functions and the corresponding class of non-differentiable ravine-type manifolds, in particular, ravine-type manifolds with corners. Such functions and manifolds are an extension of the class of smooth (differentiable) functions and manifolds, so we begin by recalling the properties of smooth objects that we need.
Below we follow to \cite{art,woll,spi,bron,rosa}.
We will consider functions on closed and open sets of the Euclidean space  
${\mathbb R}^n$.
\begin{example}
The set $[a, b]$ and the set  $(a, b)$ are the closed and the open intervals in 
${\mathbb R}$.\\
 The set $[a_1,b_1]\times \cdots \times [a_n, b_n] \subset {\mathbb R}^n $ is called a closed rectangle in ${\mathbb R}^n$, the set 
$(a_1,b_1)\times \cdots \times (a_n, b_n) \subset {\mathbb R}^n $
 is called an open rectangle.
\end{example}
    A set $U \subset {\mathbb R}^n$ is called open if for each $x \in U$ there is an open rectangle $R$ such that $x \in R \subset U$. 
 A subset $V$ of ${\mathbb R}^n$ is closed if ${\mathbb R}^n - V$ is open.

\begin{definition} 
Let $E$ and $F$ be  topological spaces. 
 Let $f$ be a one-to-one mapping of $E$ into $F$.
If both $f$ and its inverse mapping are continuous, then $f$ is called a homeomorphism, and the spaces $E$ and $F$ are called homeomorphic.
\end{definition}

\subsection{On smooth functions and smooth manifolds}
Recall some notions and results. 
\begin{definition}
Let $U $ be an open $n$-dimensional set in Euclidean space ${\mathbb R}^n$ and $f$ be a real function defined in $U $ and taking real values.
A function $f$ is called smooth if at each point in $U $ it has continuous partial derivatives of all orders, which are taken with respect to the coordinates in 
${\mathbb R}^n$. In the case $f$ is called the smooth mapping of the set $E$ into 
${\mathbb R}^n$.
\end{definition}

\begin{definition}
  Let $E$ be a set in ${\mathbb R}^n$ and let $f_1, f_2, \ldots f_n$ be $n$ smooth functions on ${\mathbb R}^n$.    The mapping $ f(x)$ of the set 
$E \subset {\mathbb R}^n$ into ${\mathbb R}^n$ at the point $(x_1, \ldots , x_n)$ has the point  $(f(x_1), \ldots , f(x_n))$  as its image of $f$.

\end{definition}

An open   $n$-dimensional ball in Euclidean space ${\mathbb R}^n$ is a set
 \begin{equation}
\label{oeb}
   B^n = \{(x_1, \ldots , x_n) | x_1^2 + \cdots  x_n^2 < 1\}.
 \end{equation} 
A set $W$ in Euclidean space ${\mathbb R}^n$ is called an open cell if it is homeomorphic to an open ball (\ref{oeb}).
\begin{definition}
An $n$-dimensional smooth manifold $M$ is a Hausdorff topological space with a covering consisting of countably many open sets $U_1, U_2, \ldots$ satisfying the following conditions:\\
   1) For each $U_i$ there is a homeomorphism $\varphi_i:U_i \to W_i$, where $W_i$ is an open cell in Euclidean space.\\
   2) If $U_i \cap U_j \neq \varnothing$, then the homeomorphism (the transition map)
$\varphi_{ij} = \varphi_j\varphi_i^{-1}$ of the set $\varphi_i(U_i \cap U_j)$ on 
$\varphi_j(U_i \cap U_j)$ is a smooth mapping.\\
The homeomorphisms ${\varphi_i}: {U_i} \to W_i$ 
are called coordinate charts, the collection $\{U_i, \varphi_i \}$ is called an atlas
of $M$.\\
The number $n$ is called the dimension ($\dim M$) of $M$ (an $m-$manifold).
\end{definition}

 \subsection{On convex functions}
  Let us recall the well-known arithmetic definition of convexity of a real function $f(x)$ of one variable on an open interval $(a, b)$.
For each two   distinct numbers $x_1, x_2$ in the interval let 
\begin{equation}
\label{simd}
   \varphi2_f(x_1, x_2) = \frac{f(x_1) - f(x_2)}{x_1 - x_2} = \varphi2_f(x_2, x_1)
 \end{equation} 
the difference quotient.

Respectively for each three   distinct numbers $x_1, x_2, x_3$ in the interval define  the quotient
\begin{equation}
\label{sim3d}
 \varphi3_f(x_1, x_2, x_3) = \frac{\varphi2_f(x_1, x_3) - \varphi2_f(x_2, x_3)}{x_1 - x_2} 
 \end{equation}.

\begin{lemma}
 The convexity of $f(x)$ is equivalent to the inequality
\begin{equation}
\label{conf}
 \varphi3_f(x_1, x_2, x_3) \ge 0
\end{equation}
for all triples of distinct numbers in the interval $(a, b)$.
\end{lemma}

   Recall the "right-handed" derivative $f^`(x_0 + 0)$ of $f(x)$ at the point $x_0$:
\begin{equation}
\label{rhsd}
 \lim_{x_1 > x_0, \;  x_1 \to x_0} \varphi2_f(x_1, x_0) = 
   \lim_{x_1 > x_0, \;  x_1 \to x_0} \frac{f(x_1) - f(x_0)}{x_1 - x_0} ,
\end{equation}
 where points $x_0, x_1, x_2$ belong to the open interval $(a, b)$.
 It is easy to see that
\[
f^\prime (x_0 + 0) \ge \varphi2_f(x_2, x_0),
\]
 and for the the "left-handed" derivative $ f^`(x_0 - 0)$
there is the inequality
\begin{equation}
\label{irlhsd}
  f^\prime (x_0 + 0) \ge f^\prime (x_0 - 0).
\end{equation}
\begin{theorem}
$f(x)$ is an convex function if, and only if, $f(x)$ has monotonically increasing one-sided derivatives. 
\end{theorem}
 \subsection{Manifolds with boundary}

Let
\begin{equation}
\label{oeb}
{\mathbb H}^m = \{x = (x_1, \ldots , x_m) \in {\mathbb R}^m| x_m \ge 0 \},
\end{equation} 
  and
\begin{equation}
 \partial{\mathbb H}^n = \{x = (x_1, \ldots , x_m) \in {\mathbb R}^m| x_m = 0 \}
\end{equation} 
 be the m-dimensional upper half space and its boundary.
\begin{definition}
A smooth $m$-manifold with boundary consists of
a  Haudorff  topological space $M$, an open cover $\{U_i\}$ 
of $M$, and a collection of homeomorphisms
\[
    \varphi_i: U_i \to \mathfrak{U}_i
\]
onto open subsets $\mathfrak{U}_i  \subset {\mathbb H}^m$, one for every $i$, such 
that, for every pair  $i, j$, the transition map
\[
\varphi_{ij} = \varphi_j\varphi_i^{-1}:  \varphi_i(U_i \cap U_j) \to \varphi_j(U_i \cap U_j)
\]
is a diffeomorphism.\\
A point $p \in M$ is called a boundary point iff it satisfies the condition: \\  
if $p \in U_i$   then $\varphi_i(p) \in  \partial{\mathbb H}^n$. \\
The set
\begin{equation}
 \partial M = \{ p \in  M | \varphi_i(p) \in \partial{\mathbb H}^n \; \forall i, \; p \in U_i \}
\end{equation} 
  of all boundary points is called the boundary of $M$.
\end{definition}

 \subsection{Manifolds with corners}

A convex combination of points $x_1, \ldots, x_n$ from ${\mathbb R}^n$ is a linear combination $\lambda_1 x_1 + \cdots \lambda_n x_n$, where 
$\lambda_1 + \cdots + \lambda_n = 1$ and $\lambda_1, \ldots , \lambda_n \ge 0$. 

\begin{theorem*} A subset $C$ of ${\mathbb R}^n$ is convex if and only if any convex combinationof points from $C$ is again in $C$.
\end{theorem*}
 Let
\begin{equation}
\label{coeb}
{\mathbb H}^{n,k} = [0, \infty)^k \times  {\mathbb R}^{n-k} =
 \{x \in {\mathbb R}^n| x_i \ge 0, 1 \le i \le k \},
\end{equation} 
Using these sets (\ref{coeb}), one can (similarly using (\ref{oeb}) for manifolds with boundary) introduce manifolds with corners.

\subsection{On one class of ravine functions}
   In Shor's monograph \cite{shor1}, a class of convex functions $f$ is considered for which in  notations of   \cite{shor1}  the inequality holds
\begin{equation}
\label{rc}
(g_f(x),x - x^*(x)) \ge \cos \varphi ||g_f(x)|| \; ||x - x^*(x)||,
\end{equation} 
Here $f$ is a convex function defined on ${\mathbb R}^n$, $g_f(x)$ is a subgradient of $f$,
$0 \le \varphi < \frac{\pi}{2}$, $x \in {\mathbb R}^n$,
 $x^*(x)$ is the point in the set oj minima of $f$ that is nearest to $x$.
If there is no $\varphi < \frac{\pi}{2}$ such that 
(\ref{rc}) holds for all $x$ in a neighborhood of the optimal point of $f$, then it is said that $f$ 
is essentially  ravine (essentially gully-shaped). 
An example of a function that  is an essentially gully-shaped   can be the example proposed by M.G. Zhurbenko \cite{shzh}: the maximum function of two positive binary quadratic forms.
Below we give an example of an infinite set of gully-shaped functions, which we construct on the basis of shifted one-dimensional Minkowski spheres.
\begin{example}
 \begin{equation}
  F_p = \max\{ f_p, g_p\}, \; p  > 1.
\end{equation}
\begin{equation}
   f_p =  |x + 1|^p + |y|^p, g_p = |x - 1|^p + |y|^p.
\end{equation}
\end{example}
 \begin{remark}. 
We will call the essentially gully-shaped functions of Example 2 Zhurbenko's essentially gully-shaped  functions (or simply Zhurbenko's ravine (gully-shaped)  functions).
\end{remark}

\section{Minkowski-Cohn surface in the Cohn parametrization.}

The ravine function described below (as well as the ravine functions described in the second part of the work) is related to the study of some problems in approximation theory, coding theory, and covering problems. Below we use  notation and results from \cite{Co:MC,GM:MM,GM:P2,GGM:PM,gopms}.

 \subsection{Minkowski balls and spheres and their points}

Let us provide definitions of the concepts used.

We consider balls of the form
\[
   D_p: \;  |x|^p + |y|^p \le 1, \; p \ge 1,
\]
and call such balls  {\it Minkowski balls}.
Continuing this, we consider the following classes of Minkowski balls and circles (one-dimensional spheres).
\begin{itemize}
   \item {\it Watson balls}: $|x|^p + |y|^p \le 1$  for $p_{0} > p \ge 2$;
  \item {\it Davis balls}: $|x|^p + |y|^p \le 1$  for $p_{0} > p \ge 2$;
  \item {\it Mordell-Chebyshev balls}: $|x|^p + |y|^p \le 1$  for $ p \ge p_{0}$;
  \end{itemize}
\begin{remark}
This classification corrects the classification given in \cite{dlmb} and more accurately reflects the historical contributions of researchers.
\end{remark}
  We consider one-dimensional  spheres of the form
\[
   {\mathbb S}_p^1 = {\mathbb S}_p : \;  |x|^p + |y|^p = 1, \; p \ge 1,
\]
and call such spheres  {\it one-dimensional Minkowski spheres} or   {\it Minkowski circles}.
\subsubsection{ Interior and exterior points of a body}
Let us recall the definition of interior and exterior points of a geometric body, starting with one-dimensional Minkowski spheres.

\begin{definition} {\bf Minkowski spheres.} Suppose a one-dimensional Minkowski sphere of the form 
${\mathbb S}_p$
is situated in the real plane, with its center at the point $a = (a_1, a_2)$ in that plane.
    Consider a point $b = (b_1, b_2)$ in the plane, where $a \neq b$. Let $L$ be the line passing through points $a$ and $b$.
Suppose the line $L$ intersects the Minkowski sphere at points $c = (c_1, c_2)$ and $d = (d_1, d_2)$. By the construction of the Minkowski sphere, the line segments connecting point $a$ to points $c$ and $d$ are equal in length. Let us denote these segments as $[a, c]$ and $[a, d]$, respectively; thus, $[a, c] = [a, d]$.
The Minkowski sphere, positioned in this way, divides the points in the plane into three sets: interior points, whose distance from the center of the sphere is less than the length of the segment $[a, c]$; exterior points, whose distance from the center of the sphere is greater than the length of the segment $[a, c]$; and points lying on the Minkowski sphere itself, whose distance from the center of the sphere is equal to the length of the segment $[a, c]$.
\end{definition}
\begin{remark}
Interior and exterior points do not belong to the Minkowski sphere.
\end{remark}

\begin{definition} {\bf Minkowski balls}.
For a Minkowski ball
\[
   D_p: \;  |x|^p + |y|^p \le 1, \; p \ge 1,
\]
whose boundary is the corresponding Minkowski sphere,
the interior and exterior points are defined in relation to this boundary sphere
\[
   {\mathbb S}_p^1 = {\mathbb S}_p : \;  |x|^p + |y|^p = 1, \; p \ge 1,
\]
A Minkowski ball consists of its interior points and the points on its boundary sphere.
\end{definition}
  \subsection{Minkowski-Cohn moduli space}
    Let Cartesian orthogonal coordinates $\tau, \sigma$ be introduced on the real plane.
Take $\tau$ and $\sigma$, $\sigma \ge \tau \ge 0$   as real slope parameters of  lattice points $(x_1, y_1), (x_2, y_2)$ in the first quadrant  on the boundary of the Minkowski sphere. The following holds:: $(x_1, y_1), (-x_2, y_2) =
 ((1 + \tau^{p})^{-\frac{1}{p}}, \tau (1 + \tau^{p})^{-\frac{1}{p}} ),
- (1 + \sigma^p)^{-\frac{1}{p}}, \sigma (1 + \sigma^p)^{-\frac{1}{p}}$.

    The Minkowski-Cohn moduli space in $\{\tau, \sigma \}$ parameterizations has the form
 \begin{equation}
\label{Delta_p(tau,sigma)}
 \Delta_p(\tau,\sigma) = (\tau + \sigma)(1 + \tau^{p})^{-\frac{1}{p}}
  (1 + \sigma^p)^{-\frac{1}{p}}, 
 \end{equation}
in the domain
 \begin{equation}
\label{M_p}
 {\mathcal M}_p: \; , 0 \le \tau \le \tau_p,  \;
 1 \leq \sigma \leq \sigma_{p}  
\end{equation}
where $\tau_p$ and $\sigma_{p}$ are defined by equations 
\begin{equation}
\label{tau_p}
   2(1 - \tau_p)^p = 1 + \tau_p^p,  \;  0 \le \tau_p < 1,
\end{equation}
 and 
\begin{equation}
\label{sigma_p}
\sigma_{p}=  (2^p - 1)^{\frac{1}{p}}
\end{equation}
of the $ \{\tau,\sigma\} $-plane, where $ \infty > p > 1$ is some real parameter.
\begin{remark}
When solving the minimization problem on the surface (\ref{Delta_p(tau,sigma)})
we have to solve minimization problem on the rectangle 
$ R_p = [0, \tau_p; 1, \sigma_{p}] $.
\end{remark}
But the minimization problem is being investigated under the condition 
\begin{equation}
\label{sigma(tau)}
   \{(1 + \tau^{p})^{-\frac{1}{p}}   - (1 + \sigma^p)^{-\frac{1}{p}}\}^p  +
 \{\tau(1 + \tau^{p})^{-\frac{1}{p}} + \sigma(1 + \sigma^p)^{-\frac{1}{p}}\}^p  = 1
\end{equation}
   Let's put 
 \begin{equation}
\label{A}
 A = (1 + \tau^{p})^{-\frac{1}{p}}   - (1 + \sigma^p)^{-\frac{1}{p}}
\end{equation}
\begin{equation}
\label{B}
B = \tau(1 + \tau^{p})^{-\frac{1}{p}} + \sigma(1 + \sigma^p)^{-\frac{1}{p}}
\end{equation}
\begin{equation}
\label{AB}
  A^p + B^p = 1
\end{equation}
From (\ref{sigma(tau)}) or from (\ref{A})-(\ref{AB}) we can define the variable $\sigma$ as a function of $\tau$.
Sometimes we will denote denote $\sigma(\tau)$ as $\sigma_{\tau}$.
\begin{lemma}
Iterative computation of $\sigma_{\tau}$ can be implemented by the formula:
\begin{equation}
\sigma_{\tau} = (1 +\sigma_{\tau}^p)^{\frac{1}{p}}((1 -A^p)^{\frac{1}{p}} - 
\tau(1 + \tau^{p})^{-\frac{1}{p}})
\end{equation}
\end{lemma}
\begin{lemma}
The level curve is determined under the condition (\ref{sigma(tau)}) by the formula
\begin{equation}
(\sigma_{\tau})_{i+1} = (1 +((\sigma_{\tau})_{i})^p)^{\frac{1}{p}}((1 -A^p)^{\frac{1}{p}} - 
\tau(1 + \tau^{p})^{-\frac{1}{p}})
\end{equation}
\end{lemma}

\bigskip

{\bf ORCID}\\
https://orcid.org/0000-0002-1586-2696 \\

\end{document}